\documentclass[12pt, a4paper]{article}

\usepackage{latexsym}
\usepackage{graphicx, float}
\usepackage{amssymb,amsfonts,amsthm,amsmath,graphicx,url}

\usepackage{hyperref}
\usepackage[labelsep=period]{caption}

\usepackage{color}

\definecolor{VeryLightBlue}{rgb}{0.9,0.9,1}
\definecolor{LightBlue}{rgb}{0.8,0.8,1}
\definecolor{MidBlue}{rgb}{0.5,0.5,1}
\definecolor{DarkBlue}{rgb}{0,0,0.6}
\definecolor{Blue}{rgb}{0,0,1}

\definecolor{Gold}{rgb}{1,0.843,0}

\definecolor{LightGreen}{rgb}{0.88,1,0.88}
\definecolor{MidGreen}{rgb}{0.6,1,0.6}
\definecolor{DarkGreen}{rgb}{0,0.6,0}

\definecolor{VeryLightYellow}{rgb}{1,1,0.9}
\definecolor{LightYellow}{rgb}{1,1,0.6}
\definecolor{MidYellow}{rgb}{1,1,0.5}
\definecolor{DarkYellow}{rgb}{1,1,0.2}

\definecolor{DarkPurple}{rgb}{.6,0,1}

\definecolor{Red}{rgb}{1,0,0}
\definecolor{VeryLightRed}{rgb}{1,0.9,0.9}
\definecolor{LightRed}{rgb}{1,0.8,0.8}
\definecolor{MidRed}{rgb}{1,0.55,0.55}

\long\def\delete#1{}

\newtheorem{theorem}{Theorem}
\newtheorem{lemma}[theorem]{Lemma}
\newtheorem{definition}{Definition}

\newcommand{\bmat}[1]{\begin{bmatrix}#1\end{bmatrix}}

\newcommand{\be}{\begin{equation}}
\newcommand{\ee}{\end{equation}}
\newcommand{\bea}{\begin{eqnarray}}
\newcommand{\eea}{\end{eqnarray}}
\newcommand{\bean}{\begin{eqnarray*}}
\newcommand{\eean}{\end{eqnarray*}}

\def\proof{\par\noindent{\textbf{Proof.}~}}

\def\qed{\hfill$\Box$\vspace{11pt}}

   \def\la{\langle}   \def\ra{\rangle}

\def\BB{{\cal B}}    \def\DD{{\cal D}}      
\def\LL{{\cal L}}        

    \def\D{{\rm D}}      \def\F{{\rm F}}

\def\bzero{{\bf 0}}          \def\bfe{{\bf e}}
        \def\bx{{\bf x}}  
      \def\bu{{\bf u}}  \def\bv{{\bf v}}  

  \def\Ga{\Gamma}       \def\Om{\Omega}  
\def\a{\alpha}   \def\b{\beta}   \def\d{\delta}   \def\g{\gamma}   \def\om{\omega}   \def\l{\lambda}   \def\s{\sigma}
\def\t{\tau}      \def\vp{\varphi}

\def\Aut{{\rm Aut}}      \def\soc{{\rm soc}}     
   \def\PSL{{\rm PSL}}   \def\PGL{{\rm PGL}}      
\def\SL{{\rm SL}}    \def\GL{{\rm GL}}        \def\AGL{{\rm AGL}}  \def\ASL{{\rm ASL}}
\def\PGL{{\rm PGL}}  \def\AG{{\rm AG}}  \def\PG{{\rm PG}}      \def\GammaL{{\rm\Gamma L}}
\def\AGammaL{{\rm A\Gamma L}}  \def\PGammaL{{\rm P\Gamma L}}    
   \def\Sp{{\rm Sp}}         
\def\F{{\rm F}}

        \def\id{{\rm id}}

\usepackage{xcolor}
\usepackage[normalem]{ulem}

\begin{document}
\openup 0.8\jot

\title{Affine flag graphs and classification of a family of symmetric graphs with complete quotients}

\author{\renewcommand{\thefootnote}{\arabic{footnote}}Yu Qing Chen\footnotemark[1] , Teng Fang\footnotemark[2] , Sanming Zhou\footnotemark[3]}

\footnotetext[1]{Department of Mathematics and Statistics, Wright State University, Dayton, OH 45435, USA}

\footnotetext[2]{School of Mathematics, Southeast University, Nanjing 211189, P. R. China}

\footnotetext[3]{{School} of Mathematics and Statistics, The University of Melbourne, Parkville, VIC 3010, Australia}

\renewcommand{\thefootnote}{}
\footnotetext{{\em E--mail addresses}: \texttt{yuqing.chen@wright.edu} (Yu Qing Chen), \texttt{tfang@seu.edu.cn} (Teng Fang), \texttt{sanming@unimelb.edu.au} (Sanming Zhou)}

\date{}

\maketitle

\begin{abstract}
A graph $\Ga$ is $G$-symmetric if $G$ is a group of automorphisms of $\Ga$ which is transitive on the set of ordered pairs of adjacent vertices of $\Ga$. If $V(\Ga)$ admits a nontrivial $G$-invariant partition $\BB$ such that for blocks $B, C \in \BB$ adjacent in the quotient graph $\Ga_{\BB}$ of $\Ga$ relative to $\BB$, exactly one vertex of $B$ has no neighbour in $C$, then $\Ga$ is called an almost multicover of $\Ga_{\BB}$. In this case an incidence structure with point set $\BB$ arises naturally, and it is a $(G, 2)$-point-transitive and $G$-block-transitive 2-design if in addition $\Ga_{\BB}$ is a complete graph. In this paper we classify all $G$-symmetric graphs $\Ga$ such that (i) $\BB$ has block size $|B| \ge 3$; (ii) $\Ga_{\BB}$ is complete and almost multi-covered by $\Ga$; (iii) the incidence structure involved is a linear space; and (iv) $G$ contains a regular normal subgroup which is elementary abelian. This classification together with earlier results in [A. Gardiner and C. E. Praeger,  Australas. J. Combin. 71 (2018) 403--426],  [M.~Giulietti et al., J. Algebraic Combin. 38  (2013) 745--765] and [T. Fang et al., Electronic J. Combin. 23 (2) (2016) P2.27] completes the classification of symmetric graphs satisfying (i) and (ii).

{\em Key words}: Symmetric graph, arc-transitive graph, linear space, flag graph
\end{abstract}

\section{Introduction}
\label{sec:int}

All graphs considered in the paper are finite and undirected. A graph is called \emph{symmetric} (or \emph{arc-transitive}) if its automorphism group is transitive on its set of arcs, where an \emph{arc} is an ordered pair of adjacent vertices. The purpose of this paper is to classify a family of symmetric graphs with complete quotients such that a certain incidence structure involved is a doubly point-transitive linear space. (A {\em linear space} \cite{Beth-Jung-Lenz} is an incidence structure of points and lines such that any point is incident with at least two lines, any line is incident with at least two points, and any two points are incident with exactly one line.) It is known that for such a linear space the group involved is either almost simple or affine. In the almost simple case the corresponding symmetric graphs have been classified in \cite{GMPZ}. In the present paper we classify the corresponding symmetric graphs in the affine case, thus completing the classification of a larger class of symmetric graphs with complete quotients.

Let $\Ga$ be a graph with vertex set $V(\Ga)$. Let $G$ be a finite group which acts on $V(\Ga)$ as a group of automorphisms of $\Ga$ (that is, $G$ preserves the adjacency and non-adjacency relations of $\Ga$). If $G$ is transitive on $V(\Ga)$ and, in its induced action, transitive on the set of arcs of $\Ga$, then $\Ga$ is said to be {\em $G$-symmetric} (or {\em $G$-arc transitive}).  If in addition $V(\Ga)$ admits a nontrivial {\em $G$-invariant partition} $\BB = \{B, C, \ldots\}$, that is, $1 < |B| < |V(\Ga)|$ and any element of $G$ maps blocks of $\BB$ to blocks of $\BB$, then $\Ga$ is called an imprimitive $G$-symmetric graph. (This occurs if and only if the stabilizer of a vertex of $\Ga$ in $G$ is not a maximal subgroup of $G$.) In this case the {\em quotient graph} of $\Ga$ relative to $\BB$, denoted by $\Ga_{\BB}$, is defined to be the graph with vertex set $\BB$ in which $B, C \in
\BB$ are adjacent if and only if there exists at least one edge of $\Ga$ between $B$ and $C$. We assume without mentioning explicitly that $\Ga_{\BB}$ has at least one edge, so that each block of $\BB$ is an independent set of $\Ga$. Denote by $\Ga(\a)$ the neighbourhood of $\a \in V(\Ga)$ in $\Ga$ and set $\Ga(B) = \cup_{\a \in
B}\Ga(\a)$. For a fixed $C \in \BB$ adjacent to $B$ in $\Ga_{\BB}$, we call
\be
\label{eq:m}
m = |\{D \in \BB: \Ga(D) \cap B = \Ga(C) \cap B\}|
\ee
the {\em multiplicity} of $\BB$. Since $\Ga$ is $G$-symmetric and $\BB$ is $G$-invariant, $|B|$, $|\Ga(C) \cap B|$ and $m$ are all independent of the choice of $B$ and $C$. If $|\Ga(C) \cap B|=|B|$ or $|\Ga(C) \cap B|=|B|-1$, then $\Ga$ is called a {\em multicover} (e.g.~\cite{LPVZ}) or {\em almost multicover} of $\Ga_{\BB}$ respectively; if in addition the edges between $B$ and $C$ form a matching, then $\Ga$ is called a {\em cover} or {\em almost cover} \cite{Zhou98} of $\Ga_{\BB}$, respectively. The reader is referred to \cite{Gardiner-Praeger95, IPZ, Li-Praeger-Zhou98, LZ, Zhou02, Zhou-EJC} for some results on imprimitive symmetric graphs, and \cite{Praeger97, Praeger00} for two excellent surveys on symmetric and highly arc-transitive graphs.

The case where $\Ga$ is an almost multicover of $\Ga_{\BB}$ is interesting because it exhibits strong connections with transitive block designs. In fact, an incidence structure, denoted by $\DD(\Ga, \BB)$, arises \cite{Zhou-EJC} naturally in this case. Its points are the blocks of $\BB$; its blocks are the images of $\BB(\a) \cup \{B\}$ under the action of $G$, where $\a \in B$ is a fixed vertex and $\BB(\a) = \{C \in \BB: \Ga(C) \cap B = B \setminus \{\a\}\}$; and its incidence relation is the set-theoretic inclusion. In general, $\DD(\Ga, \BB)$ is a 1-design of block size $m+1$ \cite[Lemma 2.2]{Zhou-EJC}. In the special case when $\Ga_{\BB}$ is a complete graph, $\DD(\Ga, \BB)$ is a $2$-$(m|B| + 1, m+1, \l)$ design with $\l = 1$ or $m+1$ admitting $G$ as a 2-point-transitive and block-transitive group of automorphisms (see \cite[Corollary 2.6]{Zhou-EJC} or \cite[Corollary 2.3]{FFXZ2}). In the case when $\l = m+1$, the corresponding graphs $\Ga$ have been classified in \cite[Theorem A]{FFXZ2}. In the case when $\l = 1$, $\DD(\Ga, \BB)$ is a $(G, 2)$-point-transitive and $G$-block-transitive linear space, and the corresponding graphs $\Ga$ have been classified in \cite[Theorem 1.1]{Gardiner-Praeger98a} (and \cite[Theorem 3.19]{Zhou-EJC} using a different approach) and \cite{GMPZ} when $\DD(\Ga, \BB)$ is trivial (that is, each line is incident with exactly two points) and nontrivial with $G$ {\em almost simple} (that is, $G$ has a nonabelian simple normal subgroup $N$ such that $N \unlhd G \le \Aut(N)$), respectively. Interesting graphs arose from such classifications, including the cross ratio graphs \cite{Gardiner-Praeger-Zhou99, Zhou-EJC} associated with finite projective lines and unitary graphs \cite{GMPZ} associated with classical Hermitian unitals. In the present paper we classify all graphs in the case when $\DD(\Ga, \BB)$ is a nontrivial linear space such that $G$ is {\em affine} (that is, $G$ contains a regular normal subgroup which is elementary abelian), thus completing the classification of all imprimitive $G$-symmetric graphs $\Ga$ with $|B| \ge 3$ such that $\Ga_{\BB}$ is  complete and almost multi-covered by $\Ga$, a project initiated in \cite{Zhou-EJC}.

Many graphs obtained from our classification are the affine flag graphs introduced in \cite{Zhou-EJC}. Let $n \ge 2$ be an integer and $q$ be a prime power. For $n \ge 3$, the Desarguesian affine space $\AG(n, q)$ is the unique $n$-dimensional affine space up to isomorphism. However, for $n = 2$, there exist other affine planes whose combinatorial parameters are the same as $\AG(2, q)$. In particular, there are four non-isomorphic affine planes of order $9$; one of them is the `exceptional nearfield affine plane', which is also called `Hughes plane'. Let $\Om(n, q)$ denote the set of point-line flags of $\AG(n, q)$. Two lines of $\AG(n, q)$ are said to be {\em intersecting} if they have a unique common point, {\em parallel} if they are on the same plane of $\AG(n, q)$ but have no point in common, and {\em skew} in the remaining case. Define $\Ga_{+}(n, q)$, $\Ga_{=}(n, q)$ and $\Ga_{\simeq}(n, q)$ to be the graphs with vertex set $\Om(n, q)$ such that two distinct flags $(\bu, L), (\bv, N) \in \Om(n,q)$ are adjacent if and only if $L$ and $N$ are intersecting, parallel and skew, respectively. These graphs were introduced in \cite[Definition 3.10]{Zhou-EJC}, where the notations $\Ga^{+}(A; n, q), \Ga^{=}(A; n, q)$ and $\Ga^{\simeq}(A; n, q)$ were used to denote them. It can be verified (see \cite[Theorem 3.14]{Zhou-EJC}) that $\Ga_{+}(n, q)$, $\Ga_{=}(n, q)$ and $\Ga_{\simeq}(n, q)$ have order $q^n (q^n - 1)/(q-1)$ and valencies $(q^n - q)(q-1)$, $q^n - q$ and $(q^{n}-q) (q^{n}-q^2)/(q-1)$, respectively. Moreover, $\Ga_{+}(n, q)$ and $\Ga_{\simeq}(n, q)$ have diameter two and girth three, while $\Ga_{=}(n, q)$ is disconnected with each component a $q^{n-1}$-partite graph with $q$ vertices in each part.

The main result in this paper is as follows.

\begin{theorem}
\label{thm:class}
Suppose that $\Ga$ is a $G$-symmetric graph admitting a nontrivial $G$-invariant
partition $\BB$ of block size $|B| \ge 3$ (where $B \in \BB$) and multiplicity $m$ such that $\Ga_{\BB}$ is a complete graph, $\Ga$ is an almost multicover of $\Ga_{\BB}$ and $\DD = \DD(\Ga, \BB)$ is a nontrivial linear space. Suppose further that $G$ contains a regular normal subgroup which is elementary abelian of order $p^d = q^n$ (with $p$ a prime and $n$ a divisor of $d$). Then one of the following occurs:
\begin{itemize}
\item[\rm (a)]
$\DD \cong \AG(n,q)$, $\SL(n, q) \unlhd G_{B}$ for some $B \in \BB$, $|B| = (q^n - 1)/(q-1)$, and $m=q-1$; moreover, the following hold:
\begin{itemize}
\item[\rm (i)]
if $n \ge 3$, then $\Ga$ is isomorphic to $\Ga_{+}(n, q)$, $\Ga_{=}(n, q)$ or $\Ga_{\simeq}(n, q)$;
\item[\rm (ii)]
if $n = 2$, then $\Ga$ is isomorphic to $\Ga_{=}(2, q)$ or belongs to a family of connected graphs with order $q^2(q+1)$;
\end{itemize}
\item[\rm (b)]
$\DD \cong \AG(2, 2)$, $G \cong \AGL(1, 4)$ or $\AGammaL(1, 4)$, $|B| = 3$, $m = 1$, and $\Ga \cong 3 \cdot K_{2,2}$ or $4 \cdot K_3$;
\item[\rm (c)]
$\DD \cong \AG(2, 4)$, $G \cong \AGammaL(1,16)$, $|B| = 5$, $m = 3$, and $\Ga \cong \Ga_{+}(2,4)$ or $\Ga_{=}(2,4)$.
\end{itemize}
\end{theorem}

The connected graphs in (a)(ii) will be defined in Definition \ref{defn:Fc}, and their valency and connectedness will be given in Lemmas \ref{lem:connected} and \ref{lem:val}, respectively.

Theorem \ref{thm:class} together with earlier results in \cite[Theorem 1.1]{Gardiner-Praeger98a} (see also \cite[Theorem 3.19]{Zhou-EJC}), \cite[Theorem 1]{GMPZ} and \cite[Theorem A]{FFXZ2} gives a complete classification of symmetric triples $(\Ga, G, \BB)$ such that $\BB$ has block size $|B| \ge 3$ and $\Ga_{\BB}$ is complete and almost multi-covered by $\Ga$.

Theorem \ref{thm:class} relies on the classification \cite{Kantor85} of doubly point-transitive linear spaces (which in turn relies on the classification of finite simple groups) and the flag graph construction introduced in \cite{Zhou-EJC}. In the proof of Theorem \ref{thm:class} we will also use a result of Cameron and Kantor \cite{CK,CK1}.

\section{Preliminaries}
\label{sec:pre}

The reader is referred to \cite{Dixon-Mortimer} and \cite{Beth-Jung-Lenz} for undefined terminology on permutation groups and combinatorial designs, respectively.

\subsection{Flag graphs}
\label{subsec:flag}

In the proof of Theorem \ref{thm:class} we will use the flag graph construction introduced in \cite{Zhou-EJC}. We give an outline of this construction for the sake of completeness of the paper.

Let $\DD$ be a $1$-design which admits a point- and block-transitive group
$G$ of automorphisms. For two points $\s, \t$ of $\DD$ and a block $L$ incident with $\s$, denote by $G_{\s}$ the stabilizer of $\s$ in $G$, by $G_{\s, \t}$ the stabilizer of $\s, \t$ in $G$ (that is, the subgroup of $G$ fixing each of $\s$ and $\t$), and by $G_{\s, L}$ the stabilizer of the flag $(\s, L)$ (that is, the setwise stabilizer of $L$ in $G_{\s}$ when $L$ is treated as a set of points). A $G$-orbit $\Om$ on the set of flags of $\DD$ is said \cite{Zhou-EJC} to be {\em feasible} with respect to $G$ if it satisfies the following conditions, where $\Om(\s)$ denotes the set of flags of $\DD$ contained in $\Om$ with point entry $\s$:
\begin{description}
\item[] (A1) $|\Om(\s)| \geq 3$;

\item[] (A2) for distinct $(\s, L), (\s, N) \in \Om(\s)$, $L \cap N = \{\s\}$;

\item[] (A3) for $(\s, L) \in \Om$, $G_{\s, L}$ is transitive on $L \setminus \{\s\}$; and

\item[] (A4) for $(\s, L) \in \Om$ and $\t \in L \setminus \{\s\}$, $G_{\s, \t}$ is transitive on $\Om(\s) \setminus \{(\s, L)\}$.
\end{description}
Since $\DD$ is $G$-point-transitive and $\Om$ is a $G$-orbit on the flags of $\DD$, the validity of these conditions is independent of the choice of point $\s$ and flag $(\s, L)$.

An ordered pair $((\s, L), (\t, N))$ of flags in $\Om$, or the corresponding $G$-orbital $\Psi = ((\s, L), (\t, N))^G$ of $\Om$, is said to be {\em compatible} \cite{Zhou-EJC} with $\Om$ if
\begin{description}
\item[] (A5) $\s \not \in N$, $\t \not \in L$ but $\s \in N'$, $\t \in L'$ for some $(\s, L'), (\t, N') \in \Om$.
\end{description}
This concept is well-defined since whenever (A5) is satisfied by some $((\s, L), (\t, N)) \in \Psi$ it is also satisfied by all other members of $\Psi$. By (A2), $(\s, L')$ and $(\t, N')$ are uniquely determined by $((\s, L), (\t, N))$.

If $\Psi$ is compatible with $\Om$ and is also {\em self-paired} (that is, $((\s, L), (\t, N)) \in \Psi$ if and only if $((\t, N), (\s, L)) \in \Psi$), then the {\em $G$-flag graph} \cite{Zhou-EJC} of $\DD$ with respect to $(\Om, \Psi)$, denoted by $\Ga(\DD, \Om, \Psi)$, is the graph with vertex set $\Om$ and arc set $\Psi$. It is proved in \cite[Theorem 1.1]{Zhou-EJC} that, for an imprimitive $G$-symmetric graph $(\Ga, \BB)$ such that $\BB$ has block size $|B| \ge 3$ (where $B \in \BB$), $\Ga$ is an almost multicover of $\Ga_{\BB}$ if and only if $\Ga$ is isomorphic to $\Ga(\DD, \Om, \Psi)$ for a $G$-point-transitive and $G$-block-transitive $1$-design $\DD$. Moreover, in this case the block size of $\DD$ is equal to $m+1$ and $\Ga_{\BB}$ has valency $m |B| $ \cite[Lemma 2.1(a)]{Zhou-EJC}, where $m$ is the multiplicity of $\BB$ defined in \eqref{eq:m}. In particular, we have the following result which is a restatement of \cite[Corollary 2.6]{Zhou-EJC}.

\begin{lemma}
\label{lem:flag graph}
(\cite[Corollary 2.6]{Zhou-EJC}) Let $s \geq 3$ be an integer and $G$ a finite group. Then the following statements are equivalent.
\begin{itemize}
\item[\rm (a)] $\Ga$ is a $G$-symmetric graph admitting a nontrivial $G$-invariant
partition $\BB$ of block size $s$ such that $\Ga_{\BB}$ is a complete graph and $\Ga$ is an almost multicover of $\Ga_{\BB}$.

\item[\rm (b)] $\Ga$ is isomorphic to $\Ga(\DD,\Om,\Psi)$ for a $(G, 2)$-point-transitive and $G$-block-transitive $2$-$(v, k, \l)$ design $\DD$ with $(v - 1)/(k - 1) = s$, a feasible $G$-orbit $\Om$ on the set of flags of $\DD$, and a self-paired $G$-orbital $\Psi$ of $\Om$ compatible with $\Om$.
\end{itemize}
Moreover, $v$ is equal to the number of vertices of the complete graph $\Ga_{\BB}$, $k - 1$ is equal to the multiplicity of $\BB$, and $G$ is faithful on the vertex set of $\Ga$ if and only if it is faithful on the point set of $\DD$.
\end{lemma}

The proof of \cite[Theorem 1.1]{Zhou-EJC} implies that for a given symmetric triple $(\Ga, G, \BB)$ the design $\DD$ in (b) is isomorphic to $\DD(\Ga, \BB)$ defined in \S\ref{sec:int}.

In the proof of Theorem \ref{thm:class} we will use the following result \cite[Lemma 2.9]{FFXZ2}: If $\DD$ is a $(G, 2)$-point-transitive and $G$-block-transitive $2$-design with point set $V$ and block size at least $3$, then a $G$-orbit $\Om=(\s, L)^G$ on the flag set of $\DD$ is feasible if and only if (i) $|\Om(\s)| \geq 3$, (ii) $L \setminus \{\s\}$ is an imprimitive block for the action of $G_\s$ on $V \setminus \{\s\}$, and (iii) $G_{\s, L}$ is transitive on $V \setminus L$.

\subsection{Some lemmas}
\label{subsec:prelimi}

\begin{lemma}
\label{lem:feasible}
Let $\DD$ be a $(G, 2)$-point-transitive linear space. Let $\Om$ be the set of flags of $\DD$.
\begin{itemize}
\item[\rm (a)] $\Om$ is the only possible feasible $G$-orbit on the set of flags of $\DD$; moreover, $\Om$ is feasible with respect to $G$ if and only if $|\Om(\s)| \geq 3$ and for some (and hence all) pair of distinct points $\s, \t$ of $\DD$, $G_{\s, \t}$ is transitive on the set of lines incident with $\s$ but not $\t$.
\item[\rm (b)] A $G$-orbital $\Psi=((\a, L), (\b, N))^G$ of $\Om$ is compatible with $\Om$ if and only if $\a\neq\b$ and $L, N\neq L(\a, \b)$, where $L(\a, \b)$ is the unique line of $\DD$ through $\a$ and $\b$.
\end{itemize}
\end{lemma}

\proof
Since $\DD$ is a $(G, 2)$-point-transitive linear space, it is necessarily $G$-flag-transitive. Hence $\Om$ is the only possible feasible $G$-orbit on the flags of $\DD$. Moreover, conditions (A2) and (A3) in \S\ref{subsec:flag} hold for $\Om$. Hence the second statement in (a) follows.

We can see that $\Psi=((\a, L), (\b, N))^G$ is compatible with $\Om$ if and only if $\a \not\in N$ and $\b \not\in L$ as we may choose $L' = N' = L(\a, \b)$ in (A5).
\qed

In the case when $\DD$ is the affine space $\AG(n,q)$ and $G\leq\AGammaL(n,q)$, the group $G$ has the following property.
\begin{lemma}
\label{lem:doubly transitive}
Let $n\geq1$ be an integer and $q$ be a prime power.
Suppose that $G \leq\AGammaL(n,q)$ is doubly point-transitive on $\DD=\AG(n,q)$. If the set of flags of $\DD$ is feasible with respect to $G$, then $G_{\bzero}$ is $2$-transitive on the set of points of $\PG(n-1, q)$, where $\bzero$ is the zero vector of $\mathbb{F}_q^n$.
\end{lemma}

\proof
Denote by $\LL$ the set of $1$-dimensional subspaces of $\mathbb{F}_q^n$. Obviously, $G_\bzero$ is transitive on $\LL$. Let $L\in \LL$ and $\t\in L \setminus \{\bzero\}$. Then $G_{\bzero, \t}$ ($\leq G_{\bzero, L}$) is transitive on $\LL\setminus\{L\}$ by the feasibility and Lemma~\ref{lem:feasible}(a). Hence $G_{\bzero}$ is $2$-transitive on $\LL$.
\qed

In the proof of Theorem \ref{thm:class} we will use the following result of Cameron and Kantor.

\begin{lemma}
\label{lem:ck}
{\em (\cite[Theorem I]{CK}; see \cite[Theorem 2.1]{CK1} for a revised version)}~Let $n \ge 2$ be an integer and $q=p^\ell$ be a prime power, where $p$ is a prime. If $H \le \GammaL(n, q)$ and $H$ is $2$-transitive on the set of points of $\PG(n-1, q)$, then one of the following holds:
\begin{itemize}
\item[\rm (a)] $\SL(n, q) \le H$;
\item[\rm (b)] $n = 4, q = 2$, and $H$ is isomorphic to the alternating group $A_7$;
\item[\rm (c)] $n = 2, q = 4$, and $H$ is a group of order 20 or 60 inducing the Frobenius group of order 20 on the set of points of $\PG(1, 4)$.
\end{itemize}
\end{lemma}

\section{Proof of Theorem~\ref{thm:class}}
\label{sec:proof}

We write the elements of $\mathbb{F}_q^n$ as column vectors. Given $\bx \in \mathbb{F}_q^n \setminus \{{\bf 0}\}$, denote the line of $\AG(n, q)$ through ${\bf 0}$ and $\bx$ by $\la \bx \ra = \{a\bx: a \in \mathbb{F}_q\}$. A typical flag of $\AG(n, q)$ can be expressed as $(\bu, \la \bx \ra+\bu)$, where $\bu \in \mathbb{F}_q^n$ and $\bx \in \mathbb{F}_q^n \setminus \{{\bf 0}\}$. {Denote $\bfe_1 = (1,0,\ldots,0)^T$, $\ldots$, $\bfe_n = (0,0,\ldots, 1)^T \in \mathbb{F}_q^n$.} A typical element of $\AGammaL(n, q)$ is denoted by
$$
t(A, \bv, \vp): \bu \mapsto A\bu^{\vp}+ \bv,\;\, \bu \in \mathbb{F}_q^n,
$$
where $A \in \GL(n, q)$, $\bv \in \mathbb{F}_q^n$, and $\vp\in \Aut(\mathbb{F}_q)$ acts componentwise {on $\mathbb{F}_q^n$}. As usual we may identify $A \in \GL(n,q)$ with $t(A, \bzero, \id)$, where $\id$ is the identity element of $\Aut(\mathbb{F}_q)$.

\subsection{Doubly transitive linear spaces (affine case)}
\label{subsec:ls}

In the proof of Theorem \ref{thm:class} we will use the classification of doubly point-transitive linear spaces \cite{Kantor85}. The group involved is either almost simple or affine. Since the almost simple case has been dealt with in \cite{GMPZ}, we focus on the affine case as required by Theorem \ref{thm:class}.

Suppose $\DD$ is a nontrivial $(G, 2)$-point-transitive linear space and $G$ contains a regular normal subgroup which is elementary abelian. Then this subgroup has order $v = p^d$, where $d \ge 1$ and $p$ is a prime. We may identify the point set of $\DD$ with some vector space $\mathbb{F}_q^n$ over $\mathbb{F}_q$ such that $q^n=p^d$ and $G\leq\AGammaL(n,q)$. The proof of \cite[Theorem 1]{Kantor85} implies that all possibilities for $(G, \DD)$ are as follows, where $\bzero$ is the zero vector of $\mathbb{F}_q^n$.
\begin{itemize}
\item[(i)] $G \le \AGammaL(1, v)$, $\DD$ is an affine space;
\item[(ii)] $\SL(n, q) \unlhd G_{\bf 0}$, $v = q^n$, $n \ge 2$;
\item[(iii)] $\Sp(n, q) \unlhd G_{\bf 0}$, $v = q^n$, $n \ge 4$ is even;
\item[(iv)] $G_2(q)' \unlhd G_{\bf 0}$, $v = q^6$, $q$ is even;
\item[(v)] $\SL(2,3) \unlhd G_{\bf 0}$ or $\SL(2,5) \unlhd G_{\bf 0}$, $\DD = \AG(2,p)$, $v=p^2$, $p = 5, 7, 11, 19, 23, 29$ or $59$;
\item[(vi)] $G_{\bf 0}$ contains a normal extraspecial subgroup $E$ of order $2^5$, $\DD = \AG(4,3)$, $v=3^4$;
\item[(vii)] $G_{\bf 0}$ contains a normal extraspecial subgroup $E$ of order $2^5$, $\DD$ is the unique `exceptional nearfield affine plane' of order $9$ (see \cite{Foulser} and \cite[pp.33-34, 229-232]{Dembowski}), and $v=3^4$;
\item[(viii)] $\SL(2, 5) \unlhd G_{\bf 0}$, $\DD$ is $\AG(2,9)$, the `exceptional nearfield affine plane' as in (vii), or $\AG(4,3)$ as in (vi), $v=3^4$;
\item[(ix)] $G_{\bf 0} = \SL(2, 13)$, $\DD = \AG(6,3)$,  $v=3^6$;
\item[(x)] $G_{\bf 0} = \SL(2, 13)$, $\DD$ is the Hering affine plane of order $27$ having $3^6$ points and $3^3 \cdot (3^3 + 1)$ lines, with $3^3$ points in each line and $3^3 + 1$ lines through each point (see \cite{Hering69} and \cite[pp.236]{Dembowski}), and $v=3^6$;
\item[(xi)] $G_{\bf 0} = \SL(2, 13)$, $\DD$ is one of the two Hering designs \cite{Hering} of order $90$ having $3^6$ points and $81 \cdot 91$ lines, with $9$ points in each line and $91$ lines through each point, and $v=3^6$.
\end{itemize}

In the rest of this section we assume that $\DD$ and $G$ are as above. We use $\Om$ to denote the set of flags of $\DD$, and $L$ to denote the unique block of $\DD$ containing $\bzero$ and $\bfe_1\in \mathbb{F}_q^n$.

\subsection{$G \le \AGammaL(1, v)$}
\label{subsec:n=1}

By \cite[Section 4]{Kantor85}, $L$ is a subfield of $\mathbb{F}_v$. Suppose $|L|=s$ and $s^{t}=v=p^d$. By Lemma~\ref{lem:feasible}, condition (A4) in \S\ref{subsec:flag} holds if and only if $G_{0,1}$ is transitive on the set of lines of $\DD$ through $\bzero$ other than $L$. (Note that in this case the zero vector $\bzero$ is the same as the zero element $0$ of $\mathbb{F}_v$.) There are $(s^t - s)/(s-1)$ such lines and thus $|G_{0,1}| \le |\AGammaL(1,v)_{0,1}| = d$. It can be verified that $d < (s^t - s)/(s-1)$ unless $(p, t, d) = (2, 2, 2)$ or $(2, 2, 4)$. In the case when $(p, t, d) = (2, 2, 2)$, we have $v = 4$ and $\DD = \AG(2, 2)$ can be viewed as the complete graph $K_4$ of four vertices, with each block $\{\s, \t\}$ treated as an edge of $K_4$ and each flag $(\s, \{\s, \t\})$ identified with the arc $(\s, \t)$ of $K_4$. Moreover, since $G$ is doubly transitive on the point set of $\DD$, we obtain $G = \AGL(1, 4)$ or $\AGammaL(1, 4)$, and one can check that $\Om$ is indeed feasible with respect to $G$. The only self-paired $G$-orbitals of $\Om$ are $\{((\s, \t), (\s', \t')): \s, \t, \s', \t' \in V(K_4)\;\mbox{pairwise distinct}\}$ and $\{((\s, \t), (\s', \t)): \s, \t, \s' \in V(K_4)\;\mbox{pairwise distinct}\}$, and the corresponding $G$-flag graphs are isomorphic to $3 \cdot K_{2,2}$ and $4 \cdot K_3$, respectively. In the case when $(p, t, d) = (2, 2, 4)$, we have $\DD = \AG(2, 4)$ and $G = \AGammaL(1,16)$, and one can verify that $\Om$ is indeed feasible with respect to $G$. By \cite[Theorem 3.13]{Zhou-EJC} and \cite[Lemma 3.9]{Zhou-EJC}, the only $G$-flag graphs from this case are $\Ga_{+}(2,4)$ and $\Ga_{=}(2,4)$.

\subsection{$\SL(n, q) \unlhd G_{\bf 0}$, $\Sp(2n, q) \unlhd G_{\bf 0}$, or $G_2(q)' \unlhd G_{\bf 0}$ ($q$ even)}
\label{subsec:n>1,SL}

This subsection covers cases (ii), (iii) and (iv). First, we have $L\subseteq\la\bfe_1\ra$ by \cite[Section 4]{Kantor85}. If $L\neq\la\bfe_1\ra$, then $G_{\bzero, L}$ is not transitive on $\mathbb{F}_q^n\setminus L$ and thus $\Om$ is not feasible by \cite[Lemma 2.9]{FFXZ2}, contradicting our assumption. Therefore, $L=\la \bfe_1\ra$, $\DD=\AG(n,q)$ and $\Om=\Om(n,q)$. Suppose that $\Om$ is feasible. Then by Lemmas~\ref{lem:doubly transitive} and \ref{lem:ck} we have $G_\bzero\geq \SL(n,q)$, and thus $\ASL(n,q)\leq G$ (the case $n=4$, $p=2$ and $G_\bzero\cong A_7$ cannot happen since $|\Sp(4, 2)|$ is not a divisor of $|A_7|$). Let
\be
\label{eq:Fnq}
\F(n,q):=\{((\s, L), (\t, N)):(\s, L), (\t, N)\in \Om,\,\s\not\in N,\,\t\not\in L\}.
\ee
Then $\F(n,q)$ is the set of ordered pairs of flags compatible with $\Om$. We use $\Psi_+(n,q)$ ($\Psi_{=}(n,q)$, $\Psi_{\simeq}(n,q)$, respectively) to denote the set of ordered pairs $((\s, L), (\t, N))$ in $\F(n,q)$ such that $L$, $N$ are intersecting (parallel, skew, respectively). Since $\ASL(n,q)\leq G$, similar to the proof of \cite[Lemma 3.9]{Zhou-EJC}, there are exactly three self-paired $G$-orbits on $\F(n,q)$ compatible with $\Om$ when $n\geq3$, namely $\Psi_+(n,q)$, $\Psi_{=}(n,q)$ and $\Psi_{\simeq}(n,q)$. Hence any $G$-flag graph of $\DD$ is isomorphic to $\Ga_+(n,q)$, $\Ga_{=}(n,q)$ or $\Ga_{\simeq}(n,q)$ when $n\geq3$.

It remains to consider the case when $n=2$. We can see that $\Psi_{=}(2,q)$ is a self-paired $G$-orbit compatible with $\Om$, and the corresponding $G$-flag graph is $\Ga_{=}(2,q)$. However, $\Psi_+(2,q)$ may not be a $G$-orbit. Let $\Psi$ be a $G$-orbit on $\Psi_+(2,q)$. Since $\ASL(2,q)\leq G$, $\Psi$ must be of the form
\begin{equation}
\label{equ:Psi}
\Psi=((\bfe_1,\la\bfe_1\ra), (c\bfe_2, \la\bfe_2\ra))^G=((-\bfe_1,\la\bfe_1+c\bfe_2\ra+c\bfe_2), (\bzero, \la\bfe_2\ra))^G
\end{equation}
for some $c\in \mathbb{F}_q^\times$. For $\vp\in\Aut(\mathbb{F}_q)$, define
$$
A_{c,\vp}=\bmat{0 & 1/{c^{\vp}} \\ c & 0},\; B_{c,\vp}=\bmat{-1 & 1/{c^\vp}\\0 & c/{c^\vp}}.
$$
If $\Psi$ is self-paired, then there exists some $t(A,\bv,\vp)\in G$ interchanging $(\bfe_1,\la\bfe_1\ra)$ and $(c\bfe_2, \la\bfe_2\ra)$, which implies that $t(A,\bv,\vp)$ stabilizes the intersecting point of $\la\bfe_1\ra$ and $\la\bfe_2\ra$ (that is, the zero vector $\bzero$), and thus $\bv=\bzero$ and $A=A_{c,\vp}$. Therefore, $\Psi$ is self-paired if and only if there exists some $t(A_{c, \d}, \bzero, \d)\in G_\bzero$ for some $\d\in\Aut(\mathbb{F}_q)$. In particular, if $p=2$, then any $G$-orbit on $\Psi_+(2,q)$ is self-paired as $A_{c, \id}\in \SL(2,q)\leq G_\bzero$.

\begin{definition}
\label{defn:Fc}
{\em Let $c\in \mathbb{F}_q^\times$. Suppose that there exists some $\d\in\Aut(\mathbb{F}_q)$ such that $t(A_{c, \d}, \bzero, \d)\in G_\bzero$. Denote by $\Ga_{G,c}(2,q)$ the $G$-flag graph $\Ga(\DD, \Om, \Psi)$, where $\Om = \Om(2,q)$ is the set of point-line flags of $\AG(2, q)$, and $\Psi$ is as in \eqref{equ:Psi}.}
\end{definition}

\begin{lemma}
\label{lem:connected}
$\Ga_{G,c}(2,q)$ is a connected graph with $q^2(q+1)$ vertices.
\end{lemma}
\proof
We use the notation above. Denote $\Ga = \Ga_{G,c}(2,q)$. By \cite[Lemma 2.11]{FFXZ2}, it suffices to prove that the group
\begin{equation}
\label{equ:g}
J:=\la G_{\bzero,\la\bfe_2\ra},g\ra
\end{equation}
is exactly $G$, where $g\in G$ interchanges $(\bzero, \la\bfe_2\ra)$ and $(-\bfe_1,\la\bfe_1+c\bfe_2\ra+c\bfe_2)$. Since $\Psi$ is self-paired by our assumption, there exists some $t(A_{c, \d}, \bzero, \d)\in G_\bzero$, and thus we can choose $g=t(B_{c, \d}, -\bfe_1, \d)$. Obviously, $g$ does not stabilize $\la\bfe_2\ra\setminus\{\bzero\}$ (as $|\la\bfe_2\ra|=q>2$). Since $G_{\bzero,\la\bfe_2\ra}$ is transitive on $\la\bfe_2\ra\setminus\{\bzero\}$ and also transitive on $\mathbb{F}_q^2\setminus\la\bfe_2\ra$ by \cite[Lemma 2.9]{FFXZ2}, $J$ is transitive on $\mathbb{F}_q^2$.

Define
$$
D_{a, b}=\bmat{a & 0 \\ b & a^{-1}}
$$
for $a\in \mathbb{F}_q^\times$ and $b \in \mathbb{F}_q$. Then $D_{a,b} \in J$ and so $h_{a,b}:=D_{a, b}^{-1}g^{-1} D_{a, b} g \in J$. It can be verified that $h_{a,b} = t(A, \bv, \id)$, where
$$
A = \bmat{\frac{a^\d}{a}+\frac{ba^\d}{c}-\frac{b^\d}{ac^\d}-\frac{bb^\d}{cc^\d}-\frac{b}{ca^\d} & \frac{ab^\d}{cc^\d}+\frac{a}{ca^\d}-\frac{aa^\d}{c} \\ -\frac{cb^\d}{ac^\d}-\frac{bb^\d}{c^\d}-\frac{b}{a^\d} & \frac{ab^\d}{c^\d}+\frac{a}{a^\d}},\;\, \bv = \left(a^\d-\frac{b^\d}{c^\d}-1, -\frac{cb^\d}{c^\d}\right)^T.
$$

\textsf{Case 1: $p>2$.} Set $a=2\in \mathbb{F}_q^\times$ and $b=0$. Then $h_{2, 0} = t\left(\bmat{1 & (1-a^2)/{c}\\0 & 1}, \bfe_1, \id\right)$. Thus $gh_{2, 0}=t\left(\bmat{-1 & (2-a^2)/{c^\d}\\0 & c/{c^\d}}, \bzero, \d\right) \in J_\bzero\setminus G_{\bzero,\la\bfe_2\ra}$. Since by \cite[Lemma 2.10]{FFXZ2}, $G_{\bzero,\la\bfe_2\ra}$ is a maximal subgroup in $G_\bzero$, we obtain $J_\bzero=G_\bzero$ and therefore $J=G$. Hence $\Ga$ is connected.

\smallskip
\textsf{Case 2: $p=2$.} Choose $g$ in \eqref{equ:g} to be $t\left(\bmat{1 & 1/{c}\\0 & 1}, \bfe_1, \id\right)$. Setting $b=c$, we obtain $h_{a, c}=t\left(\bmat{a & (a^2+a+1)/{c}\\c & a+1}, (a,c)^T, \id\right)$. Thus $h_{a,c}^2 g=t\left(\bmat{a & (a+1)^2/{c}\\c & a}, \bzero, \id\right)$. Since $q>2$, we can choose $a\in \mathbb{F}_q^\times$ such that $a+1\neq0$ and $h_{a,c}^2 g\in J_\bzero\setminus G_{\bzero,\la\bfe_2\ra}$. Therefore, $J_\bzero=G_\bzero$, $J=G$, and $\Ga$ is connected.
\qed

First, in order to determine the valency of $\Ga = \Ga_{G,c}(2,q)$, we introduce some notations. Let $\theta: \mathbb{F}_q \rightarrow \mathbb{F}_q,\; z \mapsto z^p$ be the Frobenius map of $\mathbb{F}_q$ and $\widehat{a}: \mathbb{F}_q \rightarrow \mathbb{F}_q, x\mapsto ax$ the scalar multiplication by $a\in\mathbb{F}_q^\times$. Choose $\omega$ to be a fixed primitive element of $\mathbb{F}_q$. Then $\widehat{\omega}$ generates the multiplicative group $\GL(1, q)$ and $\GammaL(1, q)=\la \widehat{\om}, \theta\ra$. Foulser \cite{Foulser64} gives a standard generating set for each subgroup of $\GammaL(1, q)$. Set $q=p^{\ell}$.

\begin{definition}
\label{def:standard form}
{\em (\cite[Definition 4.5]{LiLimPraeger})} {\em Every subgroup $M$ of $\GammaL(1,p^{\ell})$ can be uniquely presented in the form $M=\la \widehat{\om}^t, \theta^s \widehat{\om}^e\ra$ such that the following conditions all hold:
\begin{itemize}
\item[\rm (F1)] $t>0$ and $t\mid (p^{\ell}-1)$;
\item[\rm (F2)] $s>0$ and $s\mid \ell$;
\item[\rm (F3)] $0\leq e<t$ and $t\mid e(p^{\ell}-1)/(p^s-1)$.
\end{itemize}
The presentation $M=\la \widehat{\om}^t, \theta^s \widehat{\om}^e\ra$ satisfying {(F1)}--{(F3)} is said to be in {\em standard form} with {\em standard parameters} $(t, e, s)$.}
\end{definition}

The reader is referred to \cite[Lemma 4.1]{Foulser64} for the existence and uniqueness of a standard form for each subgroup of $\GammaL(1, p^{\ell})$, and to \cite{FFXZ2, FFXZ2017, Foulser64} for information about the structure of subgroups of $\GammaL(1,p^{\ell})$.

Denote $Q_u := \bmat{1&0\\0&u},\; u \in \mathbb{F}_q^\times$. Set
$$
H:=\{t(Q_u, \bzero, \d):t(Q_u, \bzero, \d)\in G_\bzero\}.
$$
Then $G_\bzero=\SL(2,q)\rtimes H$ and $H$ is exactly the stabilizer of $\bzero$, $-\bfe_1$ and $\la\bfe_2\ra$ in $G$. Define
$$
\Lambda: H\rightarrow \GammaL(1,q),\;\, t(Q_u, \bzero, \d)\mapsto t(u,0,\d),
$$
{where $t(u,0,\d): y \mapsto uy^{\d}$ for $y \in \mathbb{F}_q$}. Then $\Lambda$ is a homomorphism and is injective. Thus $H$ is isomorphic to a subgroup of $\GammaL(1,q)$. The valency of $\Ga[\Om(-\bfe_1), \Om(\bzero)]$ is equal to the length of the orbit of $\la\bfe_1+c\bfe_2\ra+c\bfe_2$ under $H$ and thus is equal to the length $\ell_c$ of the orbit of $c$ under $\Lambda(H)$. By the definition of the $G$-flag graph $\Ga$, for each $\bx\in \mathbb{F}_q^2\setminus\la\bfe_2\ra$, $(\bzero, \la\bfe_2\ra)$ has a neighbour in $\Om(\bx)$, and $(\bzero, \la\bfe_2\ra)$ has no neighbour in $\Om(\bx)$ if $\bx\in \la\bfe_2\ra$. Hence the valency of $\Ga$ is $(q^2-q)\ell_c$.

We now determine $\ell_c$. Suppose that $\Lambda(H)=\la \widehat{\om}^t, \theta^s \widehat{\om}^e\ra$ is in standard form and $c=\om^r$, where $0<r<q$. Let $i\leq \ell$ be the smallest positive integer such that $t\mid (e+r(p^s-1))\cdot\frac{p^{si}-1}{p^s-1}$. It can be verified that $i\mid \ell$, $\Lambda(H)_c$ is cyclic of order $\ell/i$, and thus $\ell_c=|\Lambda(H)|i/{\ell}=i(p^{\ell}-1)/(ts)$ (computational details are omitted). Therefore, we have proved the following result.

\begin{lemma}
\label{lem:val}
The valency of $\Ga_{G,c}(2,q)$ is equal to $iq(q-1)^2/(ts)$, where $i$, $s$ and $t$ are as above.
\end{lemma}

\subsection{$G_\bzero\leq \GL(2,p)$, $\DD = \AG(2,p)$, $p = 5,7,11,19,23,29$ or $59$}
\label{subsec:p=5,7,..}

This subsection deals with case (v) which is covered by the case when $n=2$ and $q$ is an odd prime in \S\ref{subsec:n>1,SL}. By a similar analysis we can see that there exists a self-paired $G$-orbit on $\Psi_+(2,p)$ if and only if $\bmat{1&0\\0&-1}\in G_\bzero$, {which holds} if and only if $G_\bzero=\SL(2,p)\rtimes \la C_{\ell} \ra$, where $C_{\ell} = \bmat{1&0\\0&\ell}$ for some $\ell\in\mathbb{F}_p^\times$ of even order.

Suppose that $G_\bzero=\SL(2,p)\rtimes\la C_{\ell} \ra$, where $\ell\in\mathbb{F}_p^\times$ is of even order. Then there are exactly $(p-1)/|\ell|+1$ self-paired $G$-orbits (including $\Psi_{=}(2, p)$) on $\F(2,p)$ (see \eqref{eq:Fnq}) which are compatible with $\Om$. All corresponding $G$-flag graphs excluding $\Ga_{=}(2, p)$ are connected and pairwise isomorphic with order $p^2(p+1)$ and valency $(p^2-p)\cdot|\ell|$.

\subsection{The remaining cases}
\label{subsec:remaining cases}

In this subsection we prove that no $G$-flag graph of $\DD$ arises from cases (vi)--(xi).

In cases (vi) and (vii), $G_{\bzero}/E$ is isomorphic to a subgroup of $S_5$ (see \cite[Section 2]{Kantor85}), and thus $|G_\bzero|$ divides $32\cdot 120$. Let $(\bzero, L)$ be a flag of $\DD$ and $V$ the point set of $\DD$. When $\DD=\AG(4,3)$ there are $81\cdot40$ flags, and $G_{\bzero, L}$ cannot be transitive on $V\setminus L$ as $|V\setminus L|=78$ cannot divide $|G_{\bzero, L}|=|G|/(81\cdot40)$ which is a divisor of $96$. When $\DD$ is the `exceptional nearfield affine plane' there are $90\cdot9$ flags, and $G_{\bzero, L}$ cannot be transitive on $V\setminus L$ as $|V\setminus L|=72$ cannot divide $|G_{\bzero, L}|=|G|/(90\cdot9)$ which is a divisor of $32\cdot12$. Therefore, by \cite[Lemma 2.9]{FFXZ2} there is no feasible $G$-orbit on the flag set of $\DD$. Similarly, there is no feasible $G$-orbit on the flag set of $\DD$ in cases (ix), (x) and (xi).

Finally, we deal with case (viii). Suppose that $\Om$ is a feasible $G$-orbit. Then $|G_\bzero|=240,480$ or $960$ as shown in \cite[Section 4.10, Case 1]{FFXZ2}. Similar to cases (vi) and (vii), it can be verified that $\DD$ cannot be the `exceptional nearfield affine plane' or $\AG(4,3)$. Hence $\DD=\AG(2,9)$. Now $G_\bzero$ is a subgroup of $\GammaL(2,9)$, and by Lemmas \ref{lem:doubly transitive} and \ref{lem:ck} we obtain $\SL(2,9)\leq G_\bzero$. However, this cannot happen as $|G_\bzero|$ is a divisor of $960$ but $|\SL(2,9)|$ is not. Therefore, the set of flags of $\DD$ is not a feasible $G$-orbit.

\subsection{Proof of Theorem \ref{thm:class}}
\label{subsec:pf}

Theorem \ref{thm:class} follows from Lemma \ref{lem:flag graph} and the discussion in \S\ref{subsec:ls}--\ref{subsec:remaining cases}. Note that, by Lemma \ref{lem:flag graph} and the remarks below it, $\DD$ takes the role of $\DD(\Ga, \BB)$ and so $G_{\bf 0}$ takes the role of $G_B$ for some $B \in \BB$, leading to the statements about $G_B$ in Theorem \ref{thm:class}.

\medskip
\noindent \textbf{Acknowledgements}\quad
The authors would like to thank the anonymous referees for their comments which led to improvements of presentation. The second author was supported by China Postdoctoral Science Foundation Grant 2017M621065. The third author was supported by a Future Fellowship (FT110100629) of the Australian Research Council.


{\small

\begin{thebibliography}{99}

%\bibitem{LN}
%Symplectic groups, 18.704 Supplementary Notes, March 18, 2005. \url{http://www-math.mit.edu/~dav/%symplectic.pdf}

\bibitem{Beth-Jung-Lenz}
T.~Beth, D.~Jungnickel and H.~Lenz, {\em Design Theory},
Cambridge University Press, Cambridge, 1986.

%\bibitem{Biggs}
%N.~L.~Biggs, {\em Algebraic Graph Theory} (Second edition),
%Cambridge Mathematical Library, Cambridge University Press,
%Cambridge, 1993.

%\bibitem{Biggs-White}
%N.~L.~Biggs and A.~T.~White, {\em Permutation Groups and Combinatorial
%Structures}, London Math. Soc. Lect. Note Series 33,
%Cambridge University Press, Cambridge, 1979.

%\bibitem{Buekenhout79}
%F.~Buekenhout, Diagrams for geometries and groups, {\em J. Combinatorial Theory, Ser. A},
%{\bf 27}(1979), 121-151.

%\bibitem{BDD}
%F.~Buekenhout, A.~Delandtsheer and J.~Doyen, Finite linear spaces with flag-transitive
%groups, {\em J. Combina. Theory, Ser.~A}, {\bf 49}(1988), 268-293.

%\bibitem{BDDKLS}
%F.~Buekenhout, A.~Delandtsheer and J.~Doyen, P.~B.~Kleidman, M.~W.~Liebeck and
%J.~Saxl, Linear spaces with flag-transitive automorphism groups,
%{\em Geometriae Dedicata}, {\bf 36}(1990), 89-94.

%\bibitem{Cameron81}
%P.~J.~Cameron, Finite permutation groups and finite simple groups,
%{\em Bull. London Math. Soc.}, {\bf 13}(1981), 1-22.

%\bibitem{Conder-Praeger96}
%M.~D.~E.~Conder and Cheryl E.~Praeger, Remarks on path-transitivity in finite graphs,
%{\em Europ. J. Combinatorics}, {\bf 17}(1996), 371-378.

%\bibitem{Du-Marusic-Waller98}
%S.-F.~Du, D.~Maru\v{s}i\v{c} and A.~Waller, On $2$-arc-transitive covers of
%complete graphs, {\em J. Combin. Theory Ser. B}, {\bf 74}(1998), 276--290.

%\bibitem{Atlas}
%J.~H.~Conway, R.~T.~Curtis, S.~P.~Norton, R.~A.~Parker and R.~A.~Wilson, {\em Atlas of Finite Groups}, %Clarendon Press, Oxford, 1985.

%\bibitem{BB}
%L.~M.~Batten and A.~Beutelspacher, {\em The Theory of Finite Linear Spaces}, Cambridge University Press, Cambridge, 1993.

%\bibitem{Cameron99}
%P.~J.~Cameron, {\em Permutation Groups}, Cambridge University Press, Cambridge, 1999.

\bibitem{CK}
P.~J.~Cameron and W.~M.~Kantor, $2$-transitive and antiflag transitive collineation groups of finite projective spaces, {\em J. Algebra} {\bf 60} (1979), 384--422.

\bibitem{CK1}
P.~J.~Cameron and W.~M.~Kantor, Antiflag-transitive collineation groups revisited, incomplete version, 2002, https://www.researchgate.net/publication/245428186\_Antiflag-transitive\_collineation\_groups\_revisited.

\bibitem{Dembowski}
P.~Dembowski, {\em Finite Geometries}, Springer-Verlag, Berlin, 1968.

%\bibitem{Dickson01}
%L.~E.~Dickson, {\em Linear Groups, with an Exposition of the Galois Field Theory}, Teubner, Leipzig, 1901.

\bibitem{Dixon-Mortimer}
J.~D.~Dixon and B.~Mortimer, {\em Permutation Groups},
Springer, New York, 1996.

\bibitem{FFXZ2}
T.~Fang, X.~G.~Fang, B.~Xia and S.~Zhou, A family of symmetric graphs with complete quotients, {\em Electronic J. Combin.} {\bf 23} (2) (2016), P2.27.

\bibitem{FFXZ2017}
T.~Fang, X.~G.~Fang, B.~Xia and S.~Zhou, Vertex-imprimitive symmetric graphs with exactly one edge between any two distinct blocks, {\em J. Combin. Theory Ser. A} {\bf 152} (2017), 303--340.

\bibitem{Foulser64}
D.~A.~Foulser, The flag-transitive collineation group of the finite Desarguesian affine planes,
{\em Canad. J. Math.} {\bf 16} (1964), 443--472.

\bibitem{Foulser}
D.~A.~Foulser, Solvable flag-transitive affine groups, {\em Math. Z.} {\bf 86} (1964), 191--204.

\bibitem{Gardiner-Praeger95}
A.~Gardiner and C.~E.~Praeger, A geometrical approach to
imprimitive graphs, {\em Proc. London Math. Soc.} (3) {\bf 71} (1995), 524--546.

%\bibitem{Gardiner-Praeger98}
%A.~Gardiner and Cheryl E.~Praeger, Topological covers of complete graphs,
%{\em Math. Proc. Cambridge Phil. Soc.}, {\bf 123}(1998), 549-559.

\bibitem{Gardiner-Praeger98a}
A.~Gardiner and C.~E.~Praeger, Symmetric graphs with complete
quotients, {\em Australas. J. Combin.} {\bf 71} (2018), 403--426.

\bibitem{Gardiner-Praeger-Zhou99}
A.~Gardiner, C.~E.~Praeger and S.~Zhou, Cross ratio graphs,
{\em J. London Math. Soc.} (2) {\bf 64} (2001), 257--272.

\bibitem{GMPZ}
M.~Giulietti, S.~Marcugini, F.~Pambianco and S.~Zhou, Unitary graphs and classification of a family of symmetric graphs with complete quotients, {\em J. Algebraic Combin.} {\bf 38} (2013), 745--765.

%\bibitem{Goren}
%D.~Gorenstein, R.~Lyons and R.~Solomon, {\em The Classification of the Finite Simple Groups. No. 3}, %Mathematical Surveys and Monographs 40.3, American Mathematical Society, Providence, RI, 1998.

%\bibitem{Holton-Sheehan}
%D.~A.~Holton and J.~Sheehan, {\em The Petersen Graph}, Australian Mathematical
%Society Lecture Series 7, Cambridge University Press, Cambridge, 1993.

\bibitem{Hering69}
C.~Hering, Eine nicht-desarguessche zweifach transitive affine Ebene der Ordnung 27, {\em Abh. Math. Sem. Univ. Hamburg} {\bf 34} (1969), 203--208.

\bibitem{Hering}
C.~Hering, Two new sporadic doubly transitive linear spaces, in: {\em Finite Geometries}, Lecture Notes in Pure and Applied Mathematics {\bf 103}, Marcel Dekker, New York, 1985, pp.127--129.

%\bibitem{Huppert67}
%B.~Huppert, {\em Endliche Gruppen I}, Springer-Verlag, Berlin, 1967.

%\bibitem{HKT}
%J.~W.~P.~Hirschfeld, G.~Korchm\'{a}ros and F.~Torres, Algebraic curves over a finite field, Princeton %University Press, Princeton, NJ, 2008.

%\bibitem{HP}
%D.~R.~Hughes and F.~C.~Piper, Projective planes, Springer, New York, 1973.

%\bibitem{Lu}
%H.~L\"{u}neburg, Some remarks concerning the Ree groups of type ($G_2$), {\em J.~Algebra} {\bf 3} (1966), %256--259.

\bibitem{IPZ}
M.~A.~Iranmanesh, C.~E.~Praeger and S.~Zhou, Finite symmetric graphs with two-arc transitive quotients, {\em J. Combin. Theory Ser.~B} {\bf 94} (2005), 79--99.

\bibitem{Kantor85}
W.~M.~Kantor, Homogeneous designs and geometric lattices, {\em J.~Combin. Theory Ser.~A} {\bf 38} (1985), 66--74.

%\bibitem{Kantor93}
%W.~M.~Kantor, $2$-Transitive and flag-transitive designs, in: {\em Coding Theory, %Design Theory, Group Theory} (Burlington, VT, 1990), Wiley-Interscience, Wiley, New %York, 1993, pp.13--30.

\bibitem{LiLimPraeger}
C.~H.~Li, T.~K.~Lim and C.~E.~Praeger, Homogeneous factorisations of complete graphs with edge-transitive factors, {\em J. Algebraic Combin.} {\bf 29} (2009), 107--132.

\bibitem{Li-Praeger-Zhou98}
C.~H.~Li, C.~E.~Praeger and S.~Zhou, A class of finite
symmetric graphs with 2-arc transitive quotient, {\em Math. Proc. Cambridge
Philos. Soc.} {\bf 129} (1) (2000), 19--34.

\bibitem{LPVZ}
C.~H.~Li, C.~E.~Praeger, A.~Venkatesh and
S.~Zhou, Finite locally quasiprimitive graphs, {\em Discrete
Math.}~{\bf 246} (2002), 197--218.

%\bibitem{Li-Praeger-Venkatesh-Zhou98}
%Cai Heng Li, Cheryl E.~Praeger, Akshay Venkatesh and Sanming Zhou,
%Finite locally quasiprimitive graphs, submitted.

\bibitem{LZ}
Z.~Lu and S.~Zhou, Finite symmetric graphs with $2$-arc
transitive quotients (II), {\em J. Graph Theory} {\bf 56} (2007), no.3,
167--193.

%\bibitem{MS}
%M.~Miller and J.~\v{S}ir\'{a}\v{n}, Moore graphs and beyond: A survey of the degree/diameter problem, {\em Electronic %J. Combinatorics} 20(2) (2013), \#DS14v2.

%\bibitem{Neumann-Stoy-Thompson}
%P.~M.~Neumann, G.~A.~Stoy and E.~C.~Thompson, {\em Groups and Geometry},
%Oxford University Press, Oxford, 1994.

%\bibitem{Perkel88}
%M.~Perkel, Near-polygonal graphs, {\em Ars Combinatoria}, {\bf 26(A)}(1988), 149-170.

%\bibitem{Praeger85}
%Cheryl E.~Praeger, Imprimitive symmetric graphs, {\em Ars Combinatoria},
%{\bf 19A} (1985), 149-163.

\bibitem{Praeger97}
C.~E.~Praeger,  Finite transitive permutation groups and finite vertex transitive graphs, in: G.~Hahn and G.~Sabidussi eds., {\em Graph Symmetry} (Montreal, 1996, NATO Adv. Sci. Inst. Ser. C, Math. Phys. Sci., {\bf 497}), Kluwer Academic Publishing,  Dordrecht, 1997, pp.277--318.

\bibitem{Praeger00}
C.~E.~Praeger, Finite symmetric graphs, in: L.~W.~Beineke and
R.~J.~Wilson eds., {\em Algebraic Graph Theory}, Encyclopedia of
Mathematics and Its Applications {\bf 102}, Cambridge University
Press, Cambridge, Chapter 7, pp.179-202.

%\bibitem{Sabidussi64}
%G.~Sabidussi, Vertex-transitive graphs, {\em Monatsh. Math.}, {\bf 68}(1964), 426-438.

%\bibitem{Tutte}
%W.~T.~Tutte, A family of cubical graphs,
%{\em Proc. Cambridge Philos. Soc.} {\bf 43} (1947), 459--474.

%\bibitem{Wielandt}
%H.~Wielandt, {\em Finite Permutation Groups}, Academic Press, New York, 1964.

%\bibitem{Wilson}
%R.~A.~Wilson, {\em The Finite Simple Groups}, Graduate Texts in Mathematics 251, Springer-Verlag, London, 2009.

\bibitem{Zhou02}
S.~Zhou, Imprimitive symmetric graphs, 3-arc graphs and 1-designs,
{\em Discrete Math.} {\bf 244} (2002), 521--537.

\bibitem{Zhou98}
S.~Zhou, Almost covers of $2$-arc transitive graphs, {\em
Combinatorica} {\bf 24} (2004), 731--745. [Erratum:
{\em Combinatorica} {\bf 27} (2007), 745--746]

%\bibitem{Zhou-flag}
%S.~Zhou, Symmetric graphs and flag graphs, {\em Monatshefte
%f\"{u}r Mathematik} {\bf 139} (2003), 69--81.

\bibitem{Zhou-EJC}
S.~Zhou, Constructing a class of symmetric graphs, {\em
European J. Combin.} {\bf 23} (2002), 741--760.
\end{thebibliography}
}
\end{document}